\documentclass[leqno,12pt]{article}  %putting formula numbers on the left side
\usepackage{amsmath}
\usepackage{amssymb}

\usepackage[german,english]{babel}
\usepackage{amsthm}
\author{\textsc{Elmar Grosse-Kl\"onne}}
\date{}

\theoremstyle{plain} 
\newtheorem{satz}{Theorem}[section]  %@@
\newtheorem{kor}[satz]{Corollary}  %@@
\newtheorem{lem}[satz]{Lemma}  %@@
\newtheorem{pro}[satz]{Proposition}  %@@
\theoremstyle{remark}

\theoremstyle{definition}

\newcommand{\0}{\ensuremath{\overrightarrow{0}}}

\begin{document}

% \addtolength{\textwidth}{1.0in}
% \setlength{\hoffset}{-.5in}
%
%\addtolength{\topmargin}{-34pt}
%\addtolength{\textheight}{68pt}
%\baselineskip13,6pt
%\setlength{\parindent}{0pt}
%\setlength{\parskip}{0.8ex plus 0.2ex minus 0.2ex}
%\setlength{\oddsidemargin}{3cm}
%\frenchspacing
%\sloppy
%\pagestyle{myheadings}
%\markboth{}{}

%\textwidth14cm
%\textheight9in
%\hoffset+1in
%\voffset-0,8cm
\begin{center}{\bf Locally unitary principal series representations of ${\rm
      GL}_{d+1}(F)$}\\by Elmar Grosse-Kl\"onne\\
\vspace{.5cm}
{\it To Peter Schneider on the occasion of his 60th birthday}
\end{center}

\begin{abstract} For a local field $F$ we consider tamely ramified principal
  series representations $V$ of $G={\rm GL}_{d+1}(F)$ with coefficients in a
  finite extension $K$ of ${\mathbb Q}_p$. Let $I_0$ be a pro-$p$-Iwahori
  subgroup in $G$, let ${\mathcal H}(G,I_0)$ denote the corresponding
  pro-$p$-Iwahori Hecke algebra. If $V$ is locally unitary,
  i.e. if the ${\mathcal H}(G,I_0)$-module $V^{I_0}$ admits an integral
  structure, then such an integral
  structure can be chosen in a particularly well organized manner, in
  particular its modular reduction can be made completely explicit.
\end{abstract}

\tableofcontents

\section{Introduction}

 Let $F$ be a local non-Archimedean field with finite residue field $k_F$ of
characteristic $p>0$, let ${G}={\rm GL}_{d+1}(F)$ for some $d\in{\mathbb N}$. Let $K$ be another local field which is a finite extension of ${\mathbb Q}_p$, let ${\mathfrak
  o}$ denote its ring of integers, $\pi\in {\mathfrak
  o}$ a non-zero element in its maximal ideal and $k$ its residue field. 

The general problem of deciding whether a given smooth (or, more generally, locally
algebraic) ${G}$-representation $V$ over $K$ admits a ${G}$-invariant norm --- or equivalently: a $G$-stable free ${\mathfrak o}$-sub module containing a $K$-basis of $V$ --- is of great importance for the $p$-adic local Langlands program. It is not difficult to formulate a certain {\it necessary} condition for the existence of a ${G}$-invariant norm on $V$. This has been emphasized first by Vign\'{e}ras, see also \cite{dat1}, \cite{eme}, \cite{st}, \cite{vig1}. If $V$ is a tamely ramified smooth
principal series representation and if $d=1$ then this condition turns out to also be {\it sufficient}, see \cite{vig}. Unfortunately, if $d>1$ it is unknown if this
condition is sufficient. See however \cite{univmod} for some recent progress.

In this note we consider tamely ramified smooth
principal series representations $V$ of $G$ over $K$ for general $d\in{\mathbb N}$. More precisely, we fix a maximal split torus $T$, a Borel subgroup $P$ and a pro-$p$-Iwahori subgroup $I_0$ in $G$ fixing a chamber in the apartment corresponding to $T$. We then consider a smooth $K$-valued character $\Theta$ of $T$ which is trivial on $T\cap I_0$, view it as a character of $P$ and form the smooth induction $V={\rm Ind}_P^G\Theta$. 

Let ${\mathcal
  H}(G,I_0)$ denote the pro-$p$-Iwahori Hecke algebra with coefficients in ${\mathfrak
  o}$ corresponding to $I_0$. The $K$-subspace $V^{I_0}$ of $I_0$-invariants in $V$ is naturally a module over ${\mathcal
  H}(G,I_0)\otimes_{{\mathfrak
  o}}K$. The said necessary condition for the existence of a ${G}$-invariant norm on $V$ is now equivalent with the condition that the ${\mathcal
  H}(G,I_0)\otimes_{{\mathfrak
  o}}K$-module $V^{I_0}$ admits an integral structure, i.e. an ${\mathfrak
  o}$-free ${\mathcal
  H}(G,I_0)$-sub module $L$ containing a $K$-basis of $V^{I_0}$. One might phrase this as the condition that $V$ be locally integral, or locally unitary.

It is not difficult to directly read off from
$\Theta$ whether $V$ is locally unitary. (Besides \cite{dat1} Proposition 3.2 we mention the formulation in terms of Jacquet modules as propagated by Emerton (\cite{eme}), see also section \ref{hec2} below.) We rederive this relationship here. However, the proper purpose of this paper is to provide {\it explicit} and particularly {\it well structured} ${\mathfrak
  o}$-lattices $L_{\nabla}$ in $V^{I_0}$ as above whenever $V$ is locally unitary.

Our approach is completely elementary; for example, it does not make use of the integral Bernstein basis for ${\mathcal
  H}(G,I_0)$ (e.g. \cite{vig1}). It is merely based on the investigation of
certain ${\mathbb Z}$-valued functions $\nabla$ on the finite Weyl group $W=N(T)/T$, and thus on combinatorics of $W$. We consider the canonical $K$-basis $\{f_w\}_{w\in W}$ of $V^{I_0}$ where $f_w\in V^{I_0}$ has support $PwI_0$ and satisfies $f_w(w)=1$ (we realize $W$ as a subgroup in $G$). We then ask for functions $\nabla:W\to {\mathbb Z}$ such that $L_{\nabla}=\oplus_{w\in W}(\pi)^{\nabla(w)}f_w$ is an ${\mathfrak
  o}$-lattice as desired. We show (Theorem \ref{emetr}) that whenever $V$ is locally unitary, then $V^{I_0}$ admits an ${\mathcal
  H}(G,I_0)$-stable ${\mathfrak
  o}$-lattice of this particular shape.   

The structure of the ${\mathcal
  H}(G,I_0)_k={\mathcal
  H}(G,I_0)\otimes_{{\mathfrak
  o}}k$-modules $L_{\nabla}\otimes_{{\mathfrak
  o}}k$ so obtained is then encoded in combinatorics of the (finite) Coxeter group $W$. Approaching them abstractly we suggest the notion of an ${\mathcal
  H}(G,I_0)_k$-module {\it of $W$-type} (or: a {\it reduced standard ${\mathcal
  H}(G,I_0)_k$-module}): This is an ${\mathcal
  H}(G,I_0)_k$-module $M[\theta,\sigma,\epsilon_{\bullet}]$ with $k$-basis parametrized by $W$ and whose ${\mathcal
  H}(G,I_0)_k$-structure is characterized, by means of some explicit formulae,
through a set of data $(\theta,\sigma,\epsilon_{\bullet})$ as follows:
$\theta$ is a character of $I/I_0=(T\cap I)/(T\cap I_0)$ where $I\supset I_0$
is the corresponding Iwahori subgroup; $\sigma$ is a function $\{w\in
W\,|\,\ell(ws_d)>\ell(w)\}\to\{-1,0,1\}$ where $s_d$ is the simple reflection
corresponding to an end in the Dynkin diagram, and $\ell$ is the length
function on $W$; finally, $\epsilon_{\bullet}=\{\epsilon_w\,|\,w\in W\}$ is a
set of units in $k$. (But not any such set of data
$(\theta,\sigma,\epsilon_{\bullet})$ defines an ${\mathcal
  H}(G,I_0)_k$-module $M[\theta,\sigma,\epsilon_{\bullet}]$.)  

The explicit nature of $L_{\nabla}\otimes_{{\mathfrak
  o}}k$, and more generally of an ${\mathcal
  H}(G,I_0)_k$-module of $W$-type, is particularly well suited for computing its value under a certain functor from finite dimensional ${\mathcal
  H}(G,I_0)_k$-modules to $(\varphi,\Gamma)$-modules (if $F={\mathbb Q}_p$),
see \cite{dfun}.

We intend to generalize the results of the present paper to other reductive
groups in the future. Moreover, the relationship between ${\mathcal
  H}(G,I_0)_k$-modules of $W$-type (reduced standard ${\mathcal
  H}(G,I_0)_k$-modules) and standard ${\mathcal
  H}(G,I_0)_k$-modules should be clarified.

The outline is as follows. In section \ref{fusy} we first introduce the notion of a balanced weight of length $d+1$: a $(d+1)$-tuple of integers satisfying certain boundedness conditions which later on will turn out to precisely encode the condition (on $\Theta$) for $V$ to be locally unitary. Given such a balanced weight, we show the existence of certain functions $\nabla:W\to{\mathbb Z}$ 'integrating' it. In section \ref{hec1} we introduce $V={\rm Ind}_P^G\Theta$ and show that if a function $\nabla$ 'integrates' the 'weight' associated with $\Theta$, then $L_{\nabla}$ is an ${\mathcal
  H}(G,I_0)$-stable ${\mathfrak
  o}$-lattice as desired. In section \ref{hec2} we put the results of sections \ref{fusy} and \ref{hec1} together. In section \ref{hmowt} we introduce ${\mathcal
  H}(G,I_0)_k$-modules of $W$-type.\\

{\bf Acknowledgements:} I am very grateful to the referee for critical
comments --- they helped to significantly improve the exposition.\\

\section{Functions on symmetric groups}
\label{fusy}
For a finite subset $I$ of ${\mathbb Z}_{\ge0}$ we put$$\Delta(I)=\sum_{i\in I}{i}-\frac{|I|\cdot(|I|-1)}{2}.$$

{\bf Definition:} Let $d, r\in{\mathbb N}$. We say that a sequence of integers $(n_i)_{0\le i\le
  d}=(n_0,\ldots,n_d)$ is a balanced weight of length $d+1$ and amplitude $r$ if $\sum_{i=0}^dn_i=0$ and if for each
  subset $I\subset\{0,\ldots,d\}$ we have\begin{gather}r\Delta(I)\ge \sum_{i\in I}n_{i}\ge -r\Delta(\{0,\ldots,d\}-I).\label{hyptrick}\end{gather}

\begin{lem}\label{reverwe} If $(n_i)_{0\le i\le
  d}$ is a balanced weight of length $d+1$ and amplitude $r$, then so is $(-n_{d-i})_{0\le i\le
  d}$.  
\end{lem}

{\sc Proof:} For any $I\subset \{0,\ldots,d\}$ we compute\begin{align}\Delta(I)&=\sum_{i\in
    I}i-\frac{|I|\cdot(|I|-1)}{2}\notag\\{}&=\sum_{i=0}^di-\sum_{i\notin
    I}i-d|I|-\frac{|I|^2}{2}+\frac{(d+1)|I|+d|I|}{2}\notag\\{}&=\frac{d(d+1)}{2}-\sum_{i\notin
    I}i-d|I|-\frac{|I|^2}{2}+\frac{(d+1)|I|+d|I|}{2}\notag\\{}&=d(d+1-|I|)-\sum_{i\notin
    I}i-\frac{(d+1-|I|)(d-|I|)}{2}\notag\\{}&=\sum_{i\notin
    I}(d-i) -\frac{(d+1-|I|)(d-|I|)}{2}\notag\\{}&=\Delta(\{d-i\,|\,i\in
  \{0,\ldots,d\}-I\})\notag.\end{align}Together with the
assumption $\sum_{i=0}^dn_i=0$ this shows that the set of inequalities
(\ref{hyptrick}) for $(n_i)_{0\le i\le
  d}$ is equivalent with the same set of inequalities for $(-n_{d-i})_{0\le i\le
  d}$. Namely, given $I\subset\{0,\ldots,d\}$, the inequalities (\ref{hyptrick})
for $(n_i)_{0\le i\le
  d}$ and $I$ are equivalent with the inequalities (\ref{hyptrick})
for $(-n_{d-i})_{0\le i\le
  d}$ and $\{d-i\,|\,i\in \{0,\ldots,d\}-I\}$.\hfill$\Box$\\

\begin{lem}\label{parared} Let $(n_i)_{0\le i\le
  d}$ be a balanced weight of length $d+1$ and amplitude $r$.

(a) There is a balanced weight $(\tilde{n}_i)_{0\le i\le
  d}$ of length $d+1$ and amplitude $r$ such that $\tilde{n}_0=0$ and $0\le
n_i-\tilde{n}_i\le r$ for all $1\le i\le d$. 

(b) There is a balanced weight $({m}_i)_{0\le i\le
  d-1}$ of length $d$ and amplitude $r$ such that $0\le n_i-m_{i-1}\le r$ for
each $i=1,\ldots,d$.
\end{lem}

{\sc Proof:} We first show that (b) follows from (a). Indeed, suppose we are
given $(\tilde{n}_i)_{0\le i\le
  d}$ as in (a). Then put $m_{i-1}=\tilde{n}_i$ for $i=1,\ldots,d$. We clearly have $\sum_{i=0}^{d-1}m_i=0$. Next, let $I\subset\{0,\ldots,d-1\}$. Putting $I^+=\{i+1\,|\,i\in I\}$ and
 $I^+_0=I^+\cup\{0\}$ we then find\begin{align}r\Delta(I)&=r(\sum_{i\in I}i-\frac{|I|(|I|-1)}{2})\notag\\{}&=r(\sum_{i\in I_0^+}i-|I|-\frac{|I|(|I|-1)}{2})\notag\\{}&=r(\sum_{i\in I_0^+}i-\frac{|I_0^+|(|I_0^+|-1)}{2})\notag\\{}&=r\Delta(I_0^+)\notag\\{}&\stackrel{(i)}{\ge} \sum_{i\in I^+_0}\tilde{n}_{i}=\sum_{i\in I}{m}_{i}\notag\end{align}where (i) holds true by assumption. Similarly, we find\begin{align}-r\Delta(\{0,\ldots,d-1\}-I)&=-r(\sum_{i\in \{0,\ldots,d-1\}-I}i-\frac{(d-|I|)(d-|I|-1)}{2})\notag\\{}&=-r(\sum_{i\in \{0,\ldots,d\}-I^+}i-(d-|I|)-\frac{(d-|I|)(d-|I|-1)}{2})\notag\\{}&=-r(\sum_{i\in\{0,\ldots,d\}-I^+}i-\frac{(d+1-|I^+|)(d-|I^+|)}{2})\notag\\{}&=-r\Delta(\{0,\ldots,d\}-I^+)\notag\\{}&\stackrel{(ii)}{\le} \sum_{i\in I^+}\tilde{n}_{i}=\sum_{i\in I}{m}_{i}\notag\end{align}where (ii) holds true by assumption.

Now we prove statement (a) in three steps.

Step 1: {\it For any sequence of integers $t_1,\ldots,t_d$ satisfying\begin{gather}r|I|(d-\frac{1}{2}(|I|-1))\ge
  \sum_{i\in I}t_{i}\ge\frac{1}{2}r|I|(|I|-1)\label{hypnewtrick}\end{gather}for each subset
$I\subset\{1,\ldots,d\}$, there exists another sequence of integers $\tilde{t}_1,\ldots,\tilde{t}_d$, again
satisfying formula (\ref{hypnewtrick}) for each $I\subset\{1,\ldots,d\}$ and such that
$\sum_{i=1}^d\tilde{t}_i=\frac{1}{2}rd(d-1)$ and $0\le t_i-\tilde{t}_i\le r$
for all $1\le i\le d$.}

For a subset $I\subset\{1,\ldots,d\}$ we write $I^c=\{1,\ldots,d\}-I$. Put $$\delta=\sum_{i=1}^dt_i-\frac{1}{2}rd(d-1).$$To construct $\tilde{t}_1,\ldots,\tilde{t}_d$ as desired, we put $s^{(0)}_i=t_i$ and define inductively sequences $s^{(m)}_1,\ldots,s^{(m)}_d$ for $1\le m\le \delta$ such that $0\le t_i-s_i^{(m)}\le r$, such that $0\le s_i^{(m-1)}-s_i^{(m)}\le 1$, such that $\delta-m=\sum_{i=1}^ds^{(m)}_i-\frac{1}{2}d(d-1)$ and such that for any fixed $m$ the sequence $(s^{(m)}_i)_i$ satisfies (\ref{hypnewtrick}) for each subset
$I\subset\{1,\ldots,d\}$. Once all the $(s^{(m)}_i)_i$ are constructed we may put $\tilde{t}_i=s_i^{(\delta)}$.   

Suppose $(s^{(m)}_i)_i$ have been constructed for some $m<\delta$. Let
$I_0\subset\{1,\ldots,d\}$ be maximal such that $\sum_{i\in
  I_0}s^{(m)}_{i}=\frac{1}{2}r|I_0|(|I_0|-1)$. We have\begin{gather}s^{(m)}_{i_0}<s^{(m)}_k\quad\quad\mbox{ for each }i_0\in
  I_0\mbox{ and each }k\in I_0^c.\label{maximp}\end{gather}This follows from combining
the three formulae\begin{align}\sum_{i\in
  I_0\cup\{k\}}s^{(m)}_{i}&\ge\frac{1}{2}r|I_0\cup\{k\}|(|I_0\cup\{k\}|-1)=\frac{1}{2}r|I_0|(|I_0|-1)+r|I_0|,\notag\\\sum_{i\in
  I_0}s^{(m)}_{i}&=\frac{1}{2}r|I_0|(|I_0|-1),\notag\\\sum_{i\in
  I_0-\{i_0\}}s^{(m)}_{i}&\ge\frac{1}{2}r|I_0-\{i_0\}|(|I_0-\{i_0\}|-1)=\frac{1}{2}r|I_0|(|I_0|-1)-r(|I_0|-1)\notag\end{align}(the
first one and the last one holding by hypothesis).

{\it Claim: There is some $k\in I_0^c$ such that $s_k^{(m)}+r>t_k$.}

Suppose that, on the contrary, $s_k^{(m)}+r=t_k$ for all $k\in I_0^c$. As
$(t_i)_i$ satisfies (\ref{hypnewtrick}) we then have$$
r|I_0^c|(d-\frac{1}{2}(|I_0^c|-1))\ge\sum_{k\in I_0^c}s_k^{(m)}+r$$or
equivalently$$r|I_0^c|(d-1-\frac{1}{2}(|I_0^c|-1))\ge\sum_{k\in
  I_0^c}s_k^{(m)}.$$On the other hand, as $m<\delta$ we
find\begin{align}\sum_{k\in I_0^c}s_k^{(m)}&=(\sum_{k\in
    I_0}s_k^{(m)})-\sum_{k\in
    I_0}s_k^{(m)}\notag\\{}&>\frac{1}{2}rd(d-1)-\frac{1}{2}r|I_0|(|I_0|-1)\notag\\{}&=r\sum_{n=|I_0|}^{d-1}n\notag\\{}&=r|I_0^c|(d-1-\frac{1}{2}(|I_0^c|-1)).\notag\end{align}Taken
together this is a contradiction. The claim is proven.

We choose some $k\in I_0^c$ such that $s_k^{(m)}+r>t_k$ and put $s_k^{(m+1)}=s_k^{(m)}-1$ and $s_i^{(m+1)}=s_i^{(m)}$ for
$i\in\{1,\ldots,d\}-\{k\}$. 

{\it Claim: $(s^{(m+1)}_i)_i$ satisfies the inequality on the right hand side of
(\ref{hypnewtrick}) for each $I\subset \{1,\ldots,d\}$.} 

If $k\notin I$ this follows from the inequality on the right hand side of
(\ref{hypnewtrick}) for $I$ and $(s^{(m)}_i)_i$. Similarly, if $\sum_{i\in
  I}s^{(m)}_{i}>\frac{1}{2}r|I|(|I|-1)$ the claim is obvious. Now assume that
$k\in I$ and $\sum_{i\in
  I}s^{(m)}_{i}=\frac{1}{2}r|I|(|I|-1)$. We then find some $i_0\in
I_0$ with $i_0\notin I$, because otherwise $I_0\subset I$ and hence (since $k\in I$ but
$k\notin I_0$) even $I_0\subsetneq I$, which would contradict the maximality of $I_0$ as chosen above. Formula
(\ref{maximp}) gives $s_k^{(m+1)}\ge s_{i_0}^{(m)}$, hence the inequality on the right hand side of
(\ref{hypnewtrick}) for $(I-\{k\})\cup\{i_0\}$ and $(s^{(m)}_i)_i$ implies the inequality on the right hand side of
(\ref{hypnewtrick}) for $I$ and $(s^{(m+1)}_i)_i$.

The claim is proven. All the other properties required of $(s^{(m+1)}_i)_i$ are obvious from its construction.

Step 2: {\it The sequence $t_1,\ldots,t_d$ defined by $t_i=n_i+r(d-i)$ satisfies formula (\ref{hypnewtrick}) for each subset
$I\subset\{1,\ldots,d\}$.}

Indeed, for each
$I\subset\{1,\ldots,d\}$ the formula (\ref{hypnewtrick}) for $(t_i)_{1\le
  i\le d}$ is equivalently converted into the
formula (\ref{hyptrick}) for $(n_i)_{1\le i\le d}$ by means of the following equations:\begin{align}r|I|(d-\frac{1}{2}(|I|-1))&=r\Delta(I)+\sum_{i\in
  I}r(d-i),\notag\\\frac{1}{2}r|I|(|I|-1)&=-r\Delta(\{0,\ldots,d\}-I)
+\sum_{i\in I}r(d-i).\notag\end{align}

Step 3: {\it If for the $t_i$ as in step 2 we choose $\tilde{t}_i$ as in step 1, then the sequence $(\tilde{n}_i)_{0\le i\le
  d}$ defined by $\tilde{n}_0=0$ and $\tilde{n}_i=\tilde{t}_i-r(d-i)$ for
$1\le i\le d$ satisfies the requirements
of statement (a).}

It is clear that $\tilde{n}_0=0$ and $0\le
n_i-\tilde{n}_i\le r$ for all $1\le i\le d$, as well as
$\sum_{i=0}^d\tilde{n}_i=0$. It remains to see that $(\tilde{n}_i)_{0\le i\le
  d}$ satisfies the inequalities (\ref{hyptrick}) for any
$I\subset\{0,\ldots,d\}$. If $0\notin I$ then, using the same conversion formulae as in the proof
of step 2, this follows from the fact that
$(\tilde{t}_i)_{1\le i\le d}$ satisfies formula (\ref{hyptrick}) for each
$I\subset\{1,\ldots,d\}$. If however $0\in I$ then we use the property
$\sum_{i=0}^d\tilde{n}_i=0$: it implies that, for $(\tilde{n}_i)_{0\le i\le
  d}$, the left hand (resp. right hand) side inequality of
formula (\ref{hyptrick}) for $I$ is equivalent with the right hand (resp. left hand) side inequality of
formula (\ref{hyptrick}) for $\{0,\ldots,d\}-I$, thus holds true because the latter holds true --- as we just saw.\hfill$\Box$\\ 

Let $W$ denote the finite Coxeter group of type
$A_d$. Thus, $W$ contains a set $S_0=\{s_1,\ldots, s_d\}$ of Coxeter
generators satisfying ${\rm ord}(s_is_{i+1})=3$ for $1\le i\le d-1$ and
${\rm ord}(s_is_{j+1})=2$ for $1\le i<j\le d-1$. Put $\overline{u}=s_d\cdots
s_1$. Let $\ell:W\to{\mathbb Z}_{\ge0}$ denote the length function.

It is convenient to realize $W$ as the symmetric group of the set $\{0,\ldots,d\}$
such that $s_i=(i-1,i)$ (transposition) for $1\le i\le d$. For $w\in W$ and $1\le i\le d$ we then have\begin{gather}\ell(ws_i)>\ell(w)\quad\mbox{ if and only if }\quad w(i-1)<w(i),\label{bbhi}\end{gather}see Proposition 1.5.3 in \cite{bjbr}.

Let $W'$ denote the subgroup of $W$ generated by $s_1,\ldots, s_{d-1}$. Any
element $w$ in $W$ can be uniquely written as $w=\overline{u}^iw'$ for some
$w'\in W'$, some $0\le i\le d$. We may thus define $\mu(w)=i$;
equivalently, $\mu(w)\in\{0,\ldots,d\}$ is defined by asking $\overline{u}^{-\mu(w)}w\in W'$.

\begin{satz}\label{bana} Let $(n_i)_{0\le i\le
  d}$ be a balanced weight of length $d+1$ and amplitude $r$. There exists a function $\nabla:W\to{\mathbb Z}$ such that for all $w\in
W$ we
have\begin{gather}\nabla(w)-\nabla(w\overline{u})=-n_{\mu(w)}\label{movebe}\end{gather}and
such that for all $s\in S_0$ and $w\in
W$ with $\ell(ws)>\ell(w)$ we have\begin{gather}\nabla(w)-r\le\nabla(ws)\le\nabla(w).\label{unebed}\end{gather}
\end{satz}

{\sc Proof:} We argue by induction on $d$. The case $d=1$ is
trivial. Now assume that $d\ge2$ and that we know the result for $d-1$. By
Lemma \ref{parared} we find a balanced weight $({m}_i)_{0\le i\le
  d-1}$ of length $d$ and amplitude $r$ such that $0\le n_i-m_{i-1}\le r$ for
each $i=1,\ldots,d$. Put $\overline{u}'=s_{d-1}\cdots
  s_1$. Define $\mu':W'\to\{0,\ldots,d-1\}$ by asking that for any $w\in W'$
  the element $(\overline{u}')^{-\mu'(w)}w$ of $W'$ belongs to the subgroup generated by
  $s_1,\ldots,s_{d-2}$. By induction hypothesis there is a function
$\nabla':W'\to{\mathbb Z}$
with $$\nabla'(w)-\nabla'(w\overline{u}')=-m_{\mu'(w)}$$for
all $w\in
  W'$ and $$\nabla'(w)-r\le\nabla'(ws)\le
  \nabla'(w)$$for all $w\in W', s\in \{s_1,\ldots,s_{d-1}\}$ with
  $\ell(ws)>\ell(w)$. Writing $w\in W$ uniquely as $w=w'\overline{u}^j$ with $w'\in W'$ and $0\le j\le d$ we
define$$\nabla(w)=\nabla'(w')+\sum_{t=0}^{j-1}n_{\mu(w'\overline{u}^t)}.$$That
this function $\nabla$ satisfies condition (\ref{movebe}) for all $w\in
W$ is obvious. We now
show that it satisfies condition (\ref{unebed}) for $s=s_d$ and all $w\in
W$ with $\ell(ws_d)>\ell(w)$. Write $w=w'\overline{u}^j$ with $w'\in W'$ and
$0\le j\le d$. 

If $j=d$ then $w=w'\overline{u}^d=w's_1\cdots s_{d}$ so that
$\ell(ws_d)<\ell(w)$ (since $w'\in W'$). Thus, for $j=d$ there is nothing to
prove. 

Now assume $1\le j\le d-1$. We then have
$ws_d=w\overline{u}^{-j}s_{d-j}\overline{u}^{j}=w's_{d-j}\overline{u}^{j}$
with $w's_{d-j}\in W'$, and we claim that $\ell(ws_d)>\ell(w)$ implies
$\ell(w's_{d-j})>\ell(w')$. Indeed, $\ell(ws_d)>\ell(w)$ means $w(d-1)<w(d)$, by formula (\ref{bbhi}). As $\overline{u}^j(d)=d-j$ and
$(\overline{u}')^j(d-1)=d-1-j$ this implies $w'(d-1-j)<w'(d-j)$, hence $\ell(w's_{d-j})>\ell(w')$, again by formula (\ref{bbhi}). The claim is proven.

 Moreover, for $0\le t\le j-1$ we have
$w's_{d-j}\overline{u}^{t}=w'\overline{u}^{t}s_{d-j+t}$ with $s_{d-j+t}\in
W'$. This implies
$\mu(w's_{d-j}\overline{u}^{t})=\mu(w'\overline{u}^{t})$. Therefore the claim $\nabla(w)-r\le\nabla(ws_{d})\le
  \nabla(w)$ is reduced
to the assumption $\nabla'(w')-r\le\nabla'(w's_{d-j})\le
  \nabla'(w')$. 

Finally assume that $j=0$, i.e. $w=w'\in W'$. Then $\nabla(w)=\nabla'(w)$
and \begin{align}\nabla(ws_d)&=\nabla(w\overline{u}'\overline{u}^d)\notag\\{}&=\nabla'(w\overline{u}')+\sum_{t=0}^{d-1}n_{\mu(w\overline{u}'\overline{u}^t)}.\end{align}Here
$\nabla'(w\overline{u}')=\nabla'(w)+m_{\mu'(w)}$ by the assumption on $\nabla'$. On
the other
hand $\sum_{t=0}^{d-1}n_{\mu(w\overline{u}'\overline{u}^t)}=-n_{\mu(ws_d)}$ as
$\sum_{i=0}^dn_i=0$. Now we claim that $\mu'(w)+1=\mu(ws_d)$. Indeed, we have $w(d)=d-\mu(w)$ and hence also $ws_d(d)=d-\mu(ws_d)$ for $w\in W$. Similarly, we have $w(d-1)=d-1-\mu'(w)$ and hence also $ws_d(d)=w(d-1)=d-1-\mu'(w)$ for $w\in W'$, and the claim is proven.

 Inserting
all this transforms the assumption $0\le n_{\mu(ws_d)}-m_{\mu(ws_d)-1}\le r$
into the condition (\ref{unebed}) (for $s=s_d$). 

We have proven condition (\ref{unebed}) for $s=s_d$ and all $w\in
W$ with $\ell(ws_d)>\ell(w)$. Condition (\ref{unebed}) for all $s\in S_0$ and all $w\in
W$ with $\ell(ws)>\ell(w)$ can be checked directly as well. However,
alternatively one can argue as follows.

 In the setting of section \ref{hec1}
(and in its notations) choose an arbitrary $F$ with residue field ${\mathbb
  F}_q$ (for an arbitrary $q$), and choose $K/{\mathbb Q}_p$ and $\pi\in K$
such that our present $r$ satisfies $\pi^r=q$. We use the elements
$t_{\overline{u}^i}$ of $T$ (explicitly given by formula (\ref{refreq})) to define the character
$\Theta:T\to K^{\times}$ by asking that $\Theta(t_{\overline{u}^i})=\pi^{-n_{i-1}}$ and that
$\Theta|_{T\cap I}=\theta$ be the trivial character. (This is well defined as $T$ is the direct product of $T\cap
I$ and the free abelian group on the generators $t_{\overline{u}^i}$ for $0\le i\le d$.) The implication
(iii)$\Rightarrow$(ii) in Lemma \ref{equinab}, applied to this $\Theta$, shows
that what we have proven so far is enough.\hfill$\Box$\\

\section{Hecke lattices in principal series representations I}

\label{hec1}

Fix a prime number $p$. Let $K/{\mathbb Q}_p$ be a finite extension field, ${\mathfrak
  o}$ its ring of integers and $k$ its residue field.

Let $F$ be a non-Archimedean locally compact field, ${\cal O}_F$ its ring of integers, $p_F\in{\cal O}_F$ a fixed 
prime element and $k_F={\mathbb F}_q$ its residue field with $q=p^{{\rm
    log}_pq}\in p^{\mathbb N}$ elements.

Let $G={\rm GL}_{d+1}(F)$ for some $d\in{\mathbb N}$. Let $T$ be a maximal
split torus in $G$, let $N(T)$ be its normalizer. Let $P$ be a Borel
subgroup of $G$ containing $T$, let $N$ be its unipotent radical.

Let $X$ be the Bruhat
Tits building of ${\rm PGL}_{d+1}(F)$, let $A\subset X$ be the apartment corresponding to $T$. Let $I$ be an Iwahori
subgroup of $G$ fixing a chamber $C$ in $A$, let $I_0$ denote its maximal
pro-$p$-subgroup. The (affine) reflections
in the codimension-$1$-faces of $C$ form a set $S$ of Coxeter
generators for the affine Weyl group. We view the latter as a subgroup of the
extended affine Weyl group $N(T)/T\cap I$. There is an $s_0\in S$ such that
the image of
$S_0=S-\{s_0\}$ in the finite Weyl group $W=N(T)/T$ is the set of simple reflections.  

We find elements $u,s_d\in
N(T)$ such that $uC=C$ (equivalently, $uI=Iu$, or also $uI_0=I_0u$), such that $u^{d+1}\in\{p_F\cdot {\rm id}, p_F^{-1}\cdot {\rm
  id}\}$ and such that, setting $$s_i=u^{d-i}s_du^{i-d}\quad\quad\mbox{ for }0\le i\le
d$$the set
$\{s_1,\ldots,s_d\}$ maps bijectively to $S_0$, while $\{s_0,s_1,\ldots,s_d\}$
maps bijectively to $S$; we henceforth regard these bijections as
identifications. Let
$\overline{u}=s_d\cdots s_1\in W\subset G$. Let
$\ell:W\to{\mathbb Z}_{\ge0}$ be the length function with respect to $S_0$. 

For convenience one may realize all these data explicitly, e.g. according to the following choice: $T$ consists of the diagonal matrices, $P$ consists of the upper triangular matrices, $N$ consists of the unipotent upper triangular matrices (i.e. the elements of $P$ with all diagonal entries equal to $1$). Then $W$ can be identified with the subgroup of permutation matrices in $G$. Its Coxeter generators $s_i$ for $i=1,\ldots,d$ are the block diagonal matrices$$s_i={\rm diag}(I_{i-1},\left(\begin{array}{cc}0&1\\1&0\end{array}\right),I_{d-i})$$while $u$ is written in block form as$$u=\left(\begin{array}{cc}{}&I_d\\p_F&{}\end{array}\right).$$(Here $I_m$, for $m\ge1$, always denotes the identity matrix in ${\rm GL}_{m}$.) The Iwahori group $I$ consists of the elements of ${\rm GL}_{d+1}({\mathcal O}_F)$ mapping to upper trianguler matrices in ${\rm GL}_{d+1}(k_F)$, while $I_0$ consists of the elements of $I$ whose diagonal entries map to $1\in k_F$. 

For $s\in S_0$ let $\iota_s:{\rm GL}_2(F)\to G$ denote the corresponding
embedding. For
$a\in F^{\times}$, $b\in F$ put $$h_{s}(a)=\iota_s(\left(\begin{array}{cc}a&0\\0&a^{-1}\end{array}\right)),\quad\quad
\nu_{s}(b)=\iota_s(\left(\begin{array}{cc}1&b\\0&1\end{array}\right)),\quad\quad\delta_s=\iota_s(\left(\begin{array}{cc}-1&0\\0&1\end{array}\right)).$$ We realize $W$ as a subgroup of $G$ in such a
way that $$\iota_{s}(\left(\begin{array}{cc}0&1\\1&0\end{array}\right))=s$$ for all $s\in S_0$. Notice that ${\rm Im}(\nu_s)\subset
N$ for all $s\in S_0$.

\begin{lem} (a) For $s\in S_0$ and $a\in F^{\times}$ we
  have\begin{gather}s\nu_s(a)s=h_s(a^{-1})\nu_s(a)\delta_ss\nu_s(a^{-1}).\label{hetw}\end{gather}(b)
  For $w\in W$ and $s\in S_0$ with $\ell(ws)>\ell(w)$ and for $b\in F$ we
  have\begin{gather}w\nu_s(b)w^{-1}\in N.\label{inn}\end{gather} 
\end{lem}

{\sc Proof:} Statement (a) is a straightforward computation inside
${\rm GL}_2(F)$. For statement (b) write $s=s_i$ for some $1\le i\le d$. Then
the matrix
$w\nu_s(b)w^{-1}$ has entry $b$ at the $(w(i-1),w(i))$-spot (and coincides
with the identity matrix at all other spots). As $\ell(ws_i)>\ell(w)$ implies
$w(i-1)<w(i)$ by formula (\ref{bbhi}), this implies $w\nu_s(b)w^{-1}\in N$.\hfill$\Box$\\

Let ${\rm
  ind}_{I_0}^{G}{\bf 1}_{{\mathfrak o}}$ denote the ${\mathfrak o}$-module of
${\mathfrak o}$-valued compactly supported functions $f$ on $G$ such that $f(ig)=f(g)$ for all
$g\in G$, all $i\in I_0$. It is a $G$-representation by means of the formula
$(g'f)(g)=f(gg')$ for $g,g'\in G$. Let $${\mathcal
  H}(G,I_0)={\rm End}_{{\mathfrak
    o}[G]}({\rm ind}_{I_0}^{G}{\bf 1}_{{\mathfrak o}})^{\rm op}$$denote the corresponding
pro-$p$-Iwahori Hecke algebra with coefficients in ${\mathfrak o}$. Then ${\rm
  ind}_{I_0}^{G}{\bf 1}_{{\mathfrak o}}$ is naturally a right ${\mathcal
  H}(G,I_0)$-module. For a subset $H$ of $G$ we let $\chi_H$ denote the
characteristic function of $H$. For $g\in G$ let $T_g\in {\mathcal
  H}(G,I_0)$ denote the Hecke operator corresponding to the double coset
$I_0gI_0$. It sends $f:G\to{\mathfrak o}$
to $$T_g(f):G\longrightarrow{\mathfrak o},\quad\quad h\mapsto\sum_{x\in
  I_0\backslash G}\chi_{I_0gI_0}(hx^{-1})f(x).$$In particular we have\begin{gather}T_g(\chi_{I_0})=\chi_{I_0g}=g^{-1}\chi_{I_0}\quad\quad\mbox{ if
}gI_0=I_0g.\label{hecnor}\end{gather}

Let $R$ be an ${\mathfrak o}$-algebra, let $V$ be a representation of $G$ on
an $R$-module. The submodule of $V^{I_0}$ of $I_0$-invariants in $V$ carries a
natural (left) action by the $R$-algebra ${\mathcal
  H}(G,I_0)_R={\mathcal
  H}(G,I_0)\otimes_{\mathfrak o}R$, resulting from the natural isomorphism $V^{I_0}\cong{\rm Hom}_{R[G]}(({\rm
  ind}_{I_0}^{G}{\bf 1}_{{\mathfrak o}})\otimes_{\mathfrak o}R,V)$. Explicitly, for $g\in G$ and $v\in
V^{I_0}$ the action of $T_g$ is given as follows: If the collection $\{g_j\}_j$ in $G$
is such that $I_0gI_0=\coprod_{j}I_0g_j$, then\begin{gather}T_g(v)=\sum_jg_j^{-1}v.\label{hef}\end{gather} 

Let $\overline{T}=(I\cap T)/(I_0\cap T)=I/I_0$.

Suppose we are given a character $\Theta:T\to K^{\times}$ whose
restriction $\theta=\Theta|_{I\cap T}$ to $I\cap T$ factors through
$\overline{T}$. As $\overline{T}$ is finite, $\theta$ takes values in ${\mathfrak o}^{\times}$, hence induces a character (denoted by
the same symbol) $\theta:\overline{T}\to k^{\times}$. For any $w\in W$ it defines a
homomorphism $${\theta}(wh_s(.)w^{-1}):k_F^{\times}\to
k^{\times},\quad\quad x\mapsto {\theta}(wh_s(x)w^{-1})$$ and it makes sense to compare it with the constant homomorphism
${\bf 1}$ taking all elements of $k_F^{\times}$ to $1\in k^{\times}$. Notice
in the following that ${\theta}(wh_s(.)w^{-1})={\bf 1}$ if and only if
${\theta}(wsh_s(.)sw^{-1})={\bf 1}$. For $w\in W$ and $s\in S_0$ put
$$\kappa_{w,s}=\kappa_{w,s}(\theta)=\theta(w\delta_sw^{-1})\in\{\pm 1\}.$$

Read $\Theta$ as a character of $P$ by means of the natural projection $P\to T$ and consider the smooth principal series representation$$V={\rm
  Ind}_P^G\Theta=\{f:G\to K\mbox{ locally constant
}\,|\,f(pg)=\Theta(p)f(g)\mbox{ for }g\in G, p\in P\}$$with $G$-action
$(gf)(x)=f(xg)$. For $w\in W$ let $f_w\in V$ denote the unique $I_0$-invariant
function supported on $PwI_0$ and with $f_w(w)=1$. It follows from the decomposition $G=\coprod_{w\in W}PwI_0$ that the set $\{f_w\}_{w\in W}$ is a
$K$-basis of the ${\mathcal H}(G,I_0)_K$-module
$V^{I_0}$.

\begin{lem}\label{refelem} Let $w\in W$ and $s\in S_0$, let $a\in{\mathcal O}_F$.

(a) If $\ell(ws)>\ell(w)$ and $a\notin (p_F)$ then $ws\nu_s(a)s\notin PwI_0$.

(b) If $\ell(ws)>\ell(w)$ then $v\nu_s(a)s\notin PwI_0$ for all $v\in
W-\{ws\}$.

(c) $v\nu_s(a)s\notin PwI_0$ for all $v\in
W-\{w,ws\}$.   
\end{lem}

{\sc Proof:} We have $\nu_s({\mathcal O}_F)\subset I_0$. Therefore all statements will follow from standard properties of the
decomposition $G=\coprod_{w\in W}PwI_0$, or rather the restriction of this
decomposition to ${\rm
  GL}_{d+1}({\mathcal O}_F)$; notice that this restriction projects to the
usual Bruhat decomposition of ${\rm
  GL}_{d+1}(k_F)$.

(a) The assumption $a\notin (p_F)$, i.e. $a\in{\mathcal O}_F^{\times}$, implies that
$ws\nu_s(a)s\in wIsI$, by formula (\ref{hetw}). The assumption
$\ell(ws)>\ell(w)$ implies $wIsI\subset PwsI=PwsI_0$ by standard properties of
the Bruhat decomposition, hence $wIsI\cap PwI_0=\emptyset$.

(b) Standard properties of
the Bruhat decomposition imply $vI_0s\subset PvsI_0\cup PvI_0$, as well as $vI_0s\subset PvsI_0$ if $\ell(vs)>\ell(v)$. As $\ell(ws)>\ell(w)$
and $v\ne ws$ statement (b) follows.

(c) The same argument as for (b).\hfill$\Box$\\

\begin{lem}\label{tsauffw} Let $w\in W$ and $s\in S_0$. We
  have$$T_s(f_w)=\left\{\begin{array}{l@{\quad:\quad}l}f_{ws}&\quad \mbox{  }\ell(ws)>\ell(w)       
\\qf_{ws}&\quad \mbox{
}\ell(ws)<\ell(w)\mbox{ and }{\theta}(wh_s(.)w^{-1})\ne{\bf 1}\\
qf_{ws}+\kappa_{ws,s}(q-1)f_{w}&\quad \mbox{
}\ell(ws)<\ell(w)\mbox{ and }{\theta}(wh_s(.)w^{-1})={\bf 1}\end{array}\right.$$
\end{lem}

{\sc Proof:} We have $I_0sI_0=\coprod_{a}I_0s\nu_s(a)$ where $a$ runs through a
set of representatives for $k_F$ in ${\mathcal O}_F$. For $y\in G$ we therefore compute, using formula (\ref{hef}):
\begin{align}(T_s(f_w))(y)&=(\sum_a\nu_s(a)sf_w)(y)\notag\\{}&=\sum_af_w(y\nu_s(a)s).\end{align}

Suppose first
that $\ell(ws)>\ell(w)$. For $a\notin(p_F)$ we then have
$ws\nu_s(a)s\notin PwI_0$ by Lemma \ref{refelem}, hence $f_w(ws\nu_s(a)s)=0$. On
the other hand $f_w(ws\nu_s(0)s)=f_w(w)=1$. Together we obtain
$(T_s(f_w))(ws)=1$. For $v\in W-\{ws\}$
and any $a\in{\mathcal O}_F$ we
have $v\nu_s(a)s\notin PwI_0$ by Lemma \ref{refelem}, hence $(T_s(f_w))(v)=0$. It follows that $T_s(f_w)=f_{ws}$.

Now suppose that $\ell(ws)<\ell(w)$. Then $ws\nu_s(a)sw^{-1}\in N$ for
any $a$, by formula
(\ref{inn}), hence
$f_w(ws\nu_s(a)s)=\theta(ws\nu_s(a)sw^{-1})f_w(w){=}1$. Summing up we
get$$(T_s(f_w))(ws)=\sum_af_w(ws\nu_s(a)s)=|k_F|=q.$$To compute $(T_s(f_w))(w)$ we
first notice that $f_w(w\nu_s(0)s)=f_w(ws)=0$. On the other hand, for
$a\notin(p_F)$ we find\begin{align}f_w(w\nu_s(a)s)&=f_w(wss\nu_s(a)s)\notag\\{}&\stackrel{(i)}{=}f_w(wsh_s(a^{-1})\nu_s(a)\delta_ss\nu_s(a^{-1}))\notag\\{}&=\theta(wsh_s(a^{-1})\nu_s(a)\delta_ss
w^{-1})f_w(w\nu_s(a^{-1}))\notag\\{}&\stackrel{(ii)}{=}\theta(wsh_s(a^{-1})\delta_ss
w^{-1})\notag\\{}&=\kappa_{ws,s}\theta(wsh_s(a^{-1})sw^{-1}).\notag\end{align}Here (i) uses formula (\ref{hetw}) while (ii) uses
$f_w(w\nu_s(a^{-1}))=f_w(w)=1$ as well
as $$(wsh_s(a^{-1})\nu_s(a)\delta_ss
w^{-1})\cdot(wsh_s(a^{-1})\delta_ss
w^{-1})^{-1}=ws\nu_s(a^{-1})sw^{-1}\in N,$$formula
(\ref{inn}). Now$$\sum_{a\notin(p_F)}\theta(wsh_s(a)sw^{-1})=\left\{\begin{array}{l@{\quad:\quad}l}q-1&\quad \mbox{  }{\theta}(wh_s(.)w^{-1})={\bf
  1}       
\\ 0&\quad \mbox{
}{\theta}(wh_s(.)w^{-1})\ne{\bf
  1}\end{array}\right.$$Thus
$\sum_{a\notin(p_F)}f_w(w\nu_s(a)s)=\kappa_{ws,s}(q-1)$ if
${\theta}(wh_s(.)w^{-1})={\bf 1}$, but $\sum_{a\notin(p_F)}f_w(w\nu_s(a)s)=0$ if
${\theta}(wh_s(.)w^{-1})\ne{\bf 1}$. We have shown that $(T_s(f_w))(w)=\kappa_{ws,s}(q-1)$ if
$\theta(wh_s(.)w^{-1})={\bf 1}$, but $(T_s(f_w))(w)=0$ if
$\theta(wh_s(.)w^{-1})\ne{\bf 1}$. Finally, for $v\in W-\{w,ws\}$
and $a\in{\mathcal O}_F$ we
have $v\nu_s(a)s\notin PwI_0$ by Lemma \ref{refelem}, hence $(T_s(f_w))(v)=0$. Summing up gives the formulae for $T_s(f_w)$ in
the case $\ell(ws)<\ell(w)$.\hfill$\Box$\\

As $\overline{u}$ is the unique element in $W\subset G$ lifting the image of
$u$ in $W=N(T)/T$ we have $\overline{u}^{-1}u\in T$. For $w\in W$ we define
$$t_w=w\overline{u}^{-1}uw^{-1}\in T.$$We record the formulae$$\overline{u}^{-1}u=t_{\overline{u}^0}={\rm diag}({p_F},I_{d}),$$\begin{gather}t_{\overline{u}^i}={\rm diag}(I_{d-i+1},p_F,I_{i-1})\quad\quad\mbox{ for }1\le i\le d,\label{refreq}\end{gather}In particular we notice that $t_w=t_{ws_i}$ for $2\le
i\le d$. 

\begin{lem}\label{tttuauffw} For $w\in W$ we
  have \begin{gather}T_{u^{-1}}(f_w)=\Theta(t_{w})f_{w\overline{u}^{-1}}\quad\quad\mbox{
      and }\quad\quad T_{u}(f_w)=\Theta(t^{-1}_{w\overline{u}})f_{w\overline{u}}.\label{tui}\end{gather}For
  $w\in W$ and $t\in T\cap I$ we have \begin{gather}T_{t}(f_w)=\theta(wt^{-1}w^{-1})f_{w}.\label{tt}\end{gather}
\end{lem}

{\sc Proof:} We use formula (\ref{hecnor}) in both
cases:
First, $$(T_{u^{-1}}(f_w))(w\overline{u}^{-1})=(uf_w)(w\overline{u}^{-1})=f_w(w\overline{u}^{-1}u)=\Theta(t_{w})f_w(w)=\Theta(t_{w})$$but
$$(T_{u^{-1}}(f_w))(v)=(uf_w)(v)=f_w(vu)=\Theta(vu\overline{u}^{-1}v^{-1})f_w(v\overline{u})=0$$
for $v\in W-\{w\overline{u}^{-1}\}$, hence the first one of the formulae in
(\ref{tui}); the other one is equivalent with it (or alternatively: proven in the same way). Next,$$(T_{t}(f_w))(w)=(t^{-1}f_w)(w)=f_w(wt^{-1})=\theta(wt^{-1}w^{-1})f_{w}(w)=\theta(wt^{-1}w^{-1}),$$but
$$(T_{t}(f_w))(v)=(t^{-1}f_w)(v)=f_w(vt^{-1})=\theta(vt^{-1}v^{-1})f_{w}(v)=0$$ for
$v\in W-\{w\}$, hence formula (\ref{tt}). \hfill$\Box$\\

We assume that there is some $r\in{\mathbb N}$ and some $\pi\in {\mathfrak o}$ such that
$\pi^{r}=q$ and such that $\Theta$ takes values in
the subgroup of $K^{\times}$ generated by $\pi$ and ${\mathfrak
  o}^{\times}$. Notice that, given an arbitrary $\Theta$, this can always be
achieved after passing to a suitable finite extension of $K$. Let ${\rm
  ord}_K:K\to{\mathbb Q}$ denote the order function normalized such that ${\rm ord}_K(\pi)=1$.

Suppose we are given a function $\nabla:W\to{\mathbb Z}$. For $w\in W$ we put
$g_w=\pi^{\nabla(w)}f_w$ and consider the ${\mathfrak
  o}$-submodule $$L_{\nabla}=L_{\nabla}(\Theta)=\bigoplus_{w\in W}{\mathfrak
  o}.g_w$$of $V^{I_0}$ which is ${\mathfrak
  o}$-free with basis $\{g_w\,|\,w\in W\}$. We ask under which conditions on $\nabla$ it is
stable under the action of ${\mathcal H}(G,I_0)$ on $V^{I_0}$. Consider the formulae \begin{gather}\nabla(w)-\nabla(w\overline{u})={\rm
  ord}_K(\Theta(t_{w\overline{u}})),\label{uco}\end{gather}\begin{gather}\nabla(w)-r\le\nabla(ws)\le\nabla(w).\label{wsco}\end{gather}   

\begin{lem}\label{equinab} The following conditions (i), (ii), (iii) on $\nabla$ are equivalent:

(i) $L_{\nabla}$ is stable under the action of ${\mathcal H}(G,I_0)$ on $V^{I_0}$.

(ii) $\nabla$ satisfies formula (\ref{uco}) for any $w\in W$, and it satisfies
formula (\ref{wsco}) for any $s\in S_0$ and any $w\in W$ with
$\ell(ws)>\ell(w)$.

(iii) $\nabla$ satisfies formula (\ref{uco}) for any $w\in W$, and it satisfies
formula (\ref{wsco}) for $s=s_d$ and any $w\in W$ with
$\ell(ws_d)>\ell(w)$.
\end{lem}

{\sc Proof:} For $t\in T\cap I$ and $w\in W$ it follows from Lemma
\ref{tttuauffw}
that \begin{gather}T_{t}(g_w)=\theta(wt^{-1}w^{-1})g_{w},\label{ttgw}\end{gather}\begin{gather}T_{u^{-1}}(g_w)=\pi^{\nabla(w)-\nabla(w\overline{u}^{-1})}\Theta(t_{w})g_{w\overline{u}^{-1}},\label{tugw}\end{gather}\begin{gather}T_{u}(g_w)=\pi^{\nabla(w)-\nabla(w\overline{u})}\Theta(t^{-1}_{w\overline{u}})g_{w\overline{u}}.\label{tuugw}\end{gather}For
$w\in W$ and $s\in S_0$ it follows from Lemma \ref{tsauffw} that\begin{gather}T_s(g_w)=\left\{\begin{array}{l@{\quad:\quad}l}\pi^{\nabla(w)-\nabla(ws)}g_{ws}&\quad \mbox{  }\ell(ws)>\ell(w)       
\\ \pi^{r+\nabla(w)-\nabla(ws)}g_{ws}&\quad \mbox{
}\ell(ws)<\ell(w)\mbox{ and }{\theta}(wh_s(.)w^{-1})\ne{\bf 1}\\
\pi^{r+\nabla(w)-\nabla(ws)}g_{ws}+\kappa_{ws,s}(\pi^r-1)g_{w}&\quad \mbox{
}\ell(ws)<\ell(w)\mbox{ and }{\theta}(wh_s(.)w^{-1})={\bf
  1}\end{array}\right.\label{tsgw}\end{gather}From these formulae we
immediately deduce that condition (i) implies both condition (ii) and
condition (iii) on $\nabla$. Now it is known that ${\mathcal H}(G,I_0)$ is
generated as an ${\mathfrak
  o}$-algebra by the Hecke operators $T_{t}$ for $t\in T\cap I$ together
with $T_{u^{-1}}$, $T_{u}$ and $T_{s_d}$. Thus, to show stability of
$L_{\nabla}$ under ${\mathcal H}(G,I_0)$ it is enough to show stability of
$L_{\nabla}$ under these operators. The above formulae imply that this stability is
ensured by condition (iii). Thus (i) is implied by (iii), and a fortiori by
(ii).\hfill$\Box$\\

\section{Hecke lattices in principal series representations II}

\label{hec2}

In Lemma \ref{equinab} we saw that the (particularly nice) ${\mathcal H}(G,I_0)$ stable ${\mathfrak
  o}$-lattices $L_{\nabla}$ in the ${\mathcal H}(G,I_0)_K$-module
$V^{I_0}$ for 
$V={\rm
  Ind}_P^G\Theta$ are obtained from functions $\nabla:W\to{\mathbb Z}$
satisfying the conditions stated there. We now want to explain that the
existence of such a function $\nabla$ can be directly read off from
$\Theta$. For $0\le i\le d$ put $$n_{i}=-{\rm
    ord}_K(\Theta(t_{\overline{u}^{i+1}})).$$

\begin{kor}\label{dj} If $(n_i)_{0\le i\le
  d}$ is a balanced weight of length $d+1$ and amplitude $r$ then there exists
a function $\nabla:W\to{\mathbb Z}$ such that $L_{\nabla}$ is stable under the
action of ${\mathcal H}(G,I_0)$ on $V^{I_0}$.  
\end{kor}

{\sc Proof:} By Theorem \ref{bana} there exists a function
$\nabla:W\to{\mathbb Z}$ satisfying condition (iii) of Lemma
\ref{equinab}. Thus we may conclude with that Lemma.\hfill$\Box$\\

Thus we need to decide for which $\Theta$ the collection $(n_i)_{0\le i\le
  d}$ is a balanced weight of length $d+1$ and amplitude $r$.

We now assume that $F\subset K$. We normalize the absolute value
$|.|:K^{\times}\to{\mathbb{Q}}^{\times}\subset K^{\times}$ on $K$ (and hence
its restriction to $F$) by requiring $|p_F|=q^{-1}$. Let $\delta:T\to
F^{\times}$ denote the modulus character associated with $P$,
i.e. $\delta=\prod_{\alpha\in\Phi^+}|\alpha|$ where $\Phi^+$ is the set of
positive roots. Let $N_0=N\cap I$ and $$T_+=\{t\in T\,|\,t^{-1}N_0t\subset
N_0\}.$$The group $W$ acts on the group of characters ${\rm
  Hom}(T,K^{\times})$ through its action on $T$.

\begin{satz}\label{emetr} Suppose that for all $w\in W$ and all $t\in T^+$ we
  have\begin{gather}|((w\Theta)(w\delta^{\frac{-1}{2}})\delta^{\frac{1}{2}})(t)|\le1\label{expocri}\end{gather}and that the restriction of $\Theta$ to the center of $G$ is a unitary character. Then
  $(n_i)_{0\le i\le
  d}$ is a balanced weight of length $d+1$ and amplitude $r$, and $L_{\nabla}$ is stable under the action of ${\mathcal H}(G,I_0)$ on $V^{I_0}$. 
\end{satz}

As the center of $G$ is generated by the element $\prod_{j=0}^dt_{\overline{u}^j}=p_FI_{d+1}$ (cf. formula (\ref{refreq})) together with ${\mathcal O}_F^{\times}\cdot I_{d+1}$, the condition that the restriction of $\Theta$ to the center of $G$ be a unitary character is equivalent with the condition\begin{gather}\prod_{j=0}^d|\Theta(t_{\overline{u}^j})|=1.\label{detint}\end{gather}

{\sc Proof:} (of Theorem \ref{emetr}) Recall that, for convenience, we work with the following realization: $T$ is the group of diagonal matrices, $P$ is the
group of upper triangular matrices, $s_i$ (for $1\le i\le d$) is the
$(i-1,i)$-transposition matrix and $u=\overline{u}\cdot{\rm
  diag}(p_F,1,\ldots,1)$. Thus $T_+$ is the subgroup of $T$ generated by all
$t\in\overline{T}$ (viewed as a subgroup of $T$ by means of the Teichm\"uller
character), by the scalar diagonal matrices (the center of $G$), and by all
the matrices of the form ${\rm diag}(1,\ldots,1,p_F,\ldots,p_F)$. The modulus
character is$$\delta:T\longrightarrow F^{\times},\quad\quad{\rm
  diag}(\alpha_0,\ldots,\alpha_d)\mapsto\prod_{i=0}^d|\alpha_i|^{d-2i}.$$Write
$\Theta={\rm diag}(\Theta_0,\ldots,\Theta_d)$ with characters
$\Theta_j:F^{\times}\to K^{\times}$. Reading
$W$ as the symmetric group of the set $\{0,\ldots,d\}$, formula
(\ref{expocri}) for $t={\rm diag}(\alpha_0,\ldots,\alpha_d)$
reads\begin{gather}|\prod_{i=0}^d\Theta_{\tau(i)}(\alpha_i)|\alpha_i|^{\tau(i)-i}|\le
  1\label{necconem}\end{gather}for all permutations $\tau$ of
$\{0,\ldots,d\}$. Asking formula (\ref{necconem}) for all ${\rm diag}(\alpha_0,\ldots,\alpha_d)\in T^+$ is
certainly equivalent with asking it for all
${\rm diag}(p_F^{-1},\ldots,p_F^{-1},1\ldots,1)$ and for all
${\rm diag}(1\ldots,1,p_F,\ldots,p_F)$ (and all $\tau$). This is
equivalent with asking \begin{gather}|q|^{\Delta(I)}\le|\prod_{j\in
    I}\Theta_j(p_F)|\le|q|^{-\Delta(\{0,\ldots,d\}-I)}\label{necconnew}\end{gather} for
    all $I\subset\{0,\ldots,d\}$. Indeed, the
inequalities on the left hand side of (\ref{necconnew}) are the
inequalities
(\ref{necconem}) for the
${\rm diag}(p_F^{-1},\ldots,p_F^{-1},1\ldots,1)$ and suitable $\tau$. The
inequalities on the right hand side of (\ref{necconnew}) are the
inequalities
(\ref{necconem}) for the ${\rm diag}(1\ldots,1,p_F,\ldots,p_F)$
and suitable $\tau$. Now observe that
$\Theta_j(p_F)=\Theta(t_{\overline{u}^{d+1-j}})$ and hence
$|\Theta_j(p_F)|=|\pi^{{\rm
    ord}(\Theta(t_{\overline{u}^{d+1-j}}))}|=|\pi^{-n_{d-j}}|$ for $0\le j\le
d$. We also have $|q|=|\pi^r|$. Together
with Lemma \ref{parared} we recover formula (\ref{hyptrick}). On the other hand,
formula (\ref{detint}) is just the property $\sum_{i=0}^dn_i=0$. We thus conclude with Corollary \ref{dj}.\hfill$\Box$\\

{\bf Remarks:} (1) We (formally) put $\chi=\Theta\delta^{-\frac{1}{2}}$. Let
$\overline{P}\subset G$ denote the Borel subgroup opposite to $P$. The same
arguments as in \cite{eme} page 10 show that (at least if $\chi$ is regular) for all $w\in W$ the action of $T$ on the Jacquet module
$J_{\overline{P}}(V)$ of $V$ (formed with respect to $\overline{P}$) admits a
non-zero eigenspace with character $(w\chi)\delta^{\frac{-1}{2}}$, i.e. with
character $(w\Theta)(w\delta^{\frac{-1}{2}})\delta^{\frac{-1}{2}}$. From
\cite{eme} we then deduce that the conditions in Theorem \ref{emetr} are a
necessary criterion for the existence of an integral structure in $V$.

(2) This necessary criterion has also been obtained in \cite{dat1}. Moreover,
in loc. cit. it is shown (in a much more general context) that it implies the existence of an integral structure in the
${\mathcal H}(G,I_0)$-module $V^{I_0}$. The point of Theorem \ref{emetr} is
that it explicitly describes a particularly nice such integral structure.

(3) Consider the smooth
dual ${\rm Hom}_K(V,K)^{\rm sm}$ of $V$; it is isomorphic with ${\rm
  Ind}_{P}^G\Theta^{-1}\delta$. Our conditions (\ref{expocri}) and (\ref{detint}) for $\Theta$ are
equivalent with the same conditions for
$\Theta^{-1}\delta$.\\

{\bf Remark:} Suppose we are in the setting of Corollary \ref{dj} or Theorem \ref{emetr}. Let $H$ denote a maximal compact open subgroup of $G$ containing $I$. Abstractly, $H$ is isomorphic with ${\rm GL}_{d+1}({\mathcal O}_F)$. Let ${\mathfrak o}[H].L_{\nabla}$ denote the ${\mathfrak o}[H]$-sub module of $V$ generated by $L_{\nabla}$, let $({\mathfrak o}[H].L_{\nabla})^{I_0}$ denote its ${\mathfrak o}$-sub module of $I_0$-invariants. Then one can show (we do not give the proof here) that the inclusion map $L_{\nabla}\to({\mathfrak o}[H].L_{\nabla})^{I_0}$ is surjective (and hence bijective). On the one hand this may be helpful for deciding whether $V$ contains an integral structure, i.e. a $G$-stable free ${\mathfrak o}$-sub module containing a $K$-basis of $V$. On the other hand it implies (in fact: is equivalent with it) that the induced map$$L_{\nabla}\otimes_{{\mathfrak
  o}}k\longrightarrow ({\mathfrak o}[H].L_{\nabla})\otimes_{{\mathfrak
  o}}k$$is injective. This might be a useful observation about the ${\mathcal H}(G,I_0)_k$-module $L_{\nabla}\otimes_{{\mathfrak
  o}}k$ (which we call an ${\mathcal H}(G,I_0)_k$-module of $W$-type in section \ref{hmowt}).

\section{${\mathcal H}(G,I_0)_k$-modules of $W$-type}

\label{hmowt}

We return to the setting of section \ref{hec1}. For $w\in W$ we define $$\epsilon_w=\epsilon_w(\Theta)=\pi^{-{\rm
    ord}_K(\Theta(t_{w}))}\Theta(t_{w}).$$
Let us write $W^{s_d}=\{w\in W\,|\,\ell(ws_d)>\ell(w)\}$. For a function
$\sigma:W^{s_d}\to\{-1,0,1\}$, for $w\in W$ and $i\in\{-1,0,1\}$ we understand
the condition $\sigma(w)=i$ as a shorthand for the condition [$w\in W^{s_d}$
and $\sigma(w)=i$]. 

For $w\in W$ we write $\kappa_{w}=\kappa_{ws_d,s_d}$.

Suppose that the function $\nabla:W\to{\mathbb Z}$ satisfies the equivalent conditions of
Lemma \ref{equinab}. Define a function $\sigma:W^{s_d}\to\{-1,0,1\}$ by setting\begin{gather}\sigma(w)=\left\{\begin{array}{l@{\quad:\quad}l}1&\quad \mbox{
    }\nabla(ws_d)=\nabla(w)\\ 0&\quad \mbox{
}\nabla(w)-r<\nabla(ws_d)<\nabla(w)\\
-1&\quad \mbox{
}\nabla(w)-r=\nabla(ws_d)\end{array}\right.\label{sina}\end{gather} 
The action of ${\mathcal H}(G,I_0)$ on $L_{\nabla}$ induces an action of ${\mathcal H}(G,I_0)_k={\mathcal H}(G,I_0)\otimes_{{\mathfrak
  o}}k$ on $L_{\nabla}\otimes_{{\mathfrak
  o}}k$. The ${\mathfrak
  o}$-basis $\{g_w\,|\,w\in W\}$ of $L_{\nabla}$ induces a $k$-basis
$\{g_w\,|\,w\in W\}$ of $L_{\nabla}\otimes_{{\mathfrak
  o}}k=L_{\nabla}(\Theta)\otimes_{{\mathfrak
  o}}k$ (we use the same symbols $g_w$).

\begin{kor}\label{realwty} The action of ${\mathcal H}(G,I_0)_k$ on $L_{\nabla}\otimes_{{\mathfrak
  o}}k$ is characterized through the following formulae: For $t\in T\cap I$
and $w\in W$ we
have \begin{gather}T_{t}(g_w)=\theta(wt^{-1}w^{-1})g_{w},\label{redttgw}\end{gather}\begin{gather}T_{u^{-1}}(g_w)=\epsilon_wg_{w\overline{u}^{-1}}\quad\quad\mbox{
    and }\quad\quad T_{u}(g_w)=\epsilon^{-1}_{w\overline{u}}g_{w\overline{u}},\label{redtugw}\end{gather}\begin{gather}T_{s_d}(g_w)=\left\{\begin{array}{l@{\quad:\quad}l}g_{ws_d}&\quad
      \mbox{  }[\sigma(ws_d)=-1\mbox{ and }{\theta}(wh_{s_d}(.)w^{-1})\ne{\bf 1}]\mbox{ or }\sigma(w)=1       
\\-\kappa_{w}g_{w}&\quad \mbox{
}\sigma(ws_d)\in\{0,1\}\mbox{ and }{\theta}(wh_{s_d}(.)w^{-1})={\bf 1}\\
g_{ws_d}-\kappa_wg_{w}&\quad \mbox{
}\sigma(ws_d)=-1\mbox{ and }{\theta}(wh_{s_d}(.)w^{-1})={\bf
  1}\\0&\quad \mbox{
}\mbox{ all other cases }\end{array}\right.\label{redtsgw}\end{gather}
\end{kor} 

{\sc Proof:} Formula (\ref{redttgw}) follows from formula (\ref{ttgw}). The assumption $\nabla(w\overline{u}^{-1})-\nabla(w)={\rm
  ord}_K(\theta(t_{w}))$ implies that the formulae in (\ref{redtugw})
follow from formulae (\ref{tugw}) and (\ref{tuugw}). Finally, formula (\ref{redtsgw}) follows
from formula (\ref{tsgw}) by a case by case checking.\hfill$\Box$\\ 

Forgetting their origin from some $\Theta$ and $\nabla$, we formalize the
structure of ${\mathcal H}(G,I_0)_k$-modules met in Corollary
\ref{realwty} in an independent definition.\\

{\bf Definition:} We say that an ${\mathcal H}(G,I_0)_k$-module
$M$ is of $W$-type (or: a {\it reduced standard module}) if
it is of the following form $M=M(\theta,\sigma,\epsilon_{\bullet})$. First, a $k$-vector space basis of $M$
is the set of formal symbols $g_w$ for $w\in W$. The ${\mathcal
  H}(G,I_0)_k$-action on $M$ is characterized by a character
$\theta:\overline{T}\to k^{\times}$ (which we also read as a character of $T\cap I$ by inflation), a map $\sigma:W^{s_d}\to\{-1,0,1\}$ and a
set $\epsilon_{\bullet}=\{\epsilon_{w}\}_{w\in W}$ of units $\epsilon_w\in
k^{\times}$. Namely, for $w\in W$ we define
$\kappa_w=\kappa_w(\theta)=\theta(ws_d\delta_{s_d}s_dw^{-1})\in\{\pm 1\}$. Then it
is required that for $t\in T\cap I$ and $w\in W$ formulae (\ref{redttgw}),
(\ref{redtugw}) and (\ref{redtsgw}) hold true.\\

Conversely we may begin with a character
$\theta:\overline{T}\to k^{\times}$, a map $\sigma:W^{s_d}\to\{-1,0,1\}$ and a
set $\epsilon_{\bullet}=\{\epsilon_{w}\}_{w\in W}$ of units $\epsilon_w\in
k^{\times}$ and ask:\\

{\bf Question 1:} For which set of data $\theta$, $\sigma$, $\epsilon_{\bullet}$ do formulae (\ref{redttgw}),
(\ref{redtugw}) and (\ref{redtsgw}) define an action of ${\mathcal
  H}(G,I_0)_k$ on $\oplus_{w\in W}k.g_w$ ?\\

{\bf Question 2:} For which set of data $\theta$, $\sigma$, $\epsilon_{\bullet}$ does there
exist some ${\mathcal
  H}(G,I_0)$-module $L_{\nabla}(\Theta)$ as in Corollary \ref{realwty} such
that $L_{\nabla}(\Theta)\otimes_{{\mathfrak
  o}}k\cong M(\theta,\sigma,\epsilon_{\bullet})$ as an ${\mathcal
  H}(G,I_0)_k$-module ? \\

In question 2 we regard $\theta$ as taking values in ${\mathfrak
  o}^{\times}\subset K^{\times}$ by means of the Teichm\"uller
lifting. Clearly those $\theta$, $\sigma$, $\epsilon_{\bullet}$ asked for in
question 2 belong to those $\theta$, $\sigma$, $\epsilon_{\bullet}$ asked for in
question 1. 

We do not consider question 1 in general, but provide a criterion for a
positive answer to question 2. Suppose we are given a set of data $\theta$,
$\sigma$, $\epsilon_{\bullet}$ as above.

\begin{pro}\label{silvester} Suppose that $\epsilon_w=\epsilon_{ws_i}$ for all
  $2\le i\le d$ and that there exists a function
$\partial:W\to[-r,r]\cap{\mathbb Z}$ with the
following properties:$$\sigma(w)=\left\{\begin{array}{l@{\quad:\quad}l}1&\quad \mbox{
    }w\in W^{s_d}\mbox{ and }\partial(w)=0\\ 0&\quad \mbox{
}w\in W^{s_d}\mbox{ and }0<\partial(w)<r\\
-1&\quad \mbox{
}w\in W^{s_d}\mbox{ and }\partial(w)=r\end{array}\right.$$$$\partial(ws_d)=-\partial(w)$$\begin{gather}\partial(w\overline{u}^{d-i})+\partial(ws_i\overline{u}^{d-j})=\partial(w\overline{u}^{d-j})+\partial(ws_j\overline{u}^{d-i})\label{homebed1}\end{gather}for
$1\le i<j-1<d$,
and\begin{gather}\partial(w\overline{u}^{d-i})+\partial(ws_i\overline{u}^{d-i-1})+\partial(ws_is_{i+1}\overline{u}^{d-i})=\partial(w\overline{u}^{d-i-1})+\partial(ws_{i+1}\overline{u}^{d-i})+\partial(ws_{i+1}s_i\overline{u}^{d-i-1})\label{homebed2}\end{gather}for
$1\le i<d$. 

Then there exists an extension
$\Theta:T\to K^{\times}$ of $\theta$ and a function $\nabla:W\to{\mathbb Z}$
as before such
that we have an isomorphism of ${\mathcal
  H}(G,I_0)_k$-modules $L_{\nabla}(\Theta)\otimes_{{\mathfrak
  o}}k\cong M(\theta,\sigma,\epsilon_{\bullet})$.
\end{pro}

{\sc Proof:} {\it Step 1:} Let $w,v\in W$. Choose a (not necessarily reduced) expression
$v=s_{i_1}\cdots s_{i_r}$ (with $i_m\in\{1,\ldots,d\}$) and
put$$\partial(w,v)=\sum_{m=1}^r\partial(ws_{i_1}\cdots
s_{i_{m-1}}\overline{u}^{d-i_m}).$${\it Claim: This definition does not depend
  on the chosen expression $s_{i_1}\cdots s_{i_r}$ for $v$.}

Indeed, it follows from hypothesis (\ref{homebed1}) that for $1\le i<j-1<d$ we
have $\partial(w,s_is_j)=\partial(w,s_js_i)$ where on either side we use the
expression of $s_is_j=s_js_i$ as indicated. Similarly, it follows from
hypothesis (\ref{homebed2}) that for $1\le i<d$ we
have $\partial(w,s_is_{i+1}s_i)=\partial(w,s_{i+1}s_is_{i+1})$ where on either side we use the
expression of $s_is_{i+1}s_i=s_{i+1}s_is_{i+1}$ as indicated. Finally, for
$1\le i\le d$ we have $\partial(w,s_is_i)=0$ where we use the expression
$s_is_i$ for the element $s_is_i=s_i^2=1\in W$: this follows from the
definition of $\partial$ and from
$s_i\overline{u}^{d-i}=\overline{u}^{d-i}s_{d}$. Thus we see that our
definition of $\partial(w,v)$ (viewed as a function in $v\in W$, with fixed
$w\in W$) respects the defining relations for the Coxeter group $W$. Iterated application
implies the stated claim.  

{\it Step 2:} The definition of $\partial(w,v)$ implies
$\partial(w,v)+\partial(wv,x)=\partial(w,vx)$ for $v,w,x\in W$. Therefore
there is a function $\nabla:W\to{\mathbb Z}$, uniquely determined up to
addition of a constant function $W\to{\mathbb Z}$, such
that$$\nabla(w)-\nabla(wv)=\partial(w,v)\quad\quad\mbox{ for all }v,w\in
W.$$It has the following properties. First, it fulfils formula
(\ref{sina}). Next, we have \begin{gather}\nabla(w)-\nabla(w\overline{u})=\nabla(ws_i)-\nabla(ws_i\overline{u})\quad\quad\mbox{
  for }w\in W\mbox{ and }1\le i\le d-1.\label{hausfor0}\end{gather}\begin{gather}\nabla(w\overline{u}^{-1})-\nabla(w)=\nabla(w\overline{u}^{-1}s_i)-\nabla(ws_i)\quad\quad\mbox{
  for }w\in W\mbox{ and }2\le i\le d.\label{hausfor}\end{gather}These formulae
are equivalent, as $s_{i}\overline{u}=\overline{u}s_{i+1}$ for $1\le i\le
d-1$. To see that they hold true we compute
\begin{align}\nabla(w)-\nabla(ws_i)&=\partial(w,s_i)\notag\\{}&=\partial(w\overline{u}^{d-i})\notag\\{}&=\partial(w\overline{u},s_{i+1})\notag\\{}&=\nabla(w\overline{u})-\nabla(w\overline{u}s_{i+1})\notag\\{}&=\nabla(w\overline{u})-\nabla(ws_i\overline{u})\end{align}
and formula (\ref{hausfor0}) follows.

{\it Step 3:} For $w\in W$ we
define$$\Theta(t_w)=\pi^{\nabla(w\overline{u}^{-1})-\nabla(w)}\epsilon_w \in
K^{\times}.$$Formula (\ref{hausfor}) together with our assumption on the
$\epsilon_w$ implies that this is well defined, because for $w, w'\in W$ we
have $t_w=t_{w'}$ if and only if $w^{-1}w'$ belongs to the subgroup of $W$
generated by $s_2,\ldots,s_{d}$. As $T/T\cap I$ is
freely generated by the $t_w$ this defines a character $\Theta:T\to
K^{\times}$ extending $T\cap I\to\overline{T}\stackrel{\theta}{\to}
k^{\times}\subset K^{\times}$, as desired.\hfill$\Box$\\

\begin{kor} Assume that $d\le 2$. If we have $\epsilon_w=\epsilon_{ws_i}$ for all
  $2\le i\le d$ then there exists an extension
$\Theta:T\to K^{\times}$ of $\theta$ and a function $\nabla:W\to{\mathbb Z}$ such
that we have an isomorphism of ${\mathcal
  H}(G,I_0)_k$-modules $L_{\nabla}(\Theta)\otimes_{{\mathfrak
  o}}k\cong M(\theta,\sigma,\epsilon_{\bullet})$.
\end{kor}

{\sc Proof:} Choose a function
$\partial:W^{s_d}\to[0,r]\cap{\mathbb Z}$ such that$$\partial(w)=0\mbox{ if }\sigma(w)=1,\quad\quad 0<\partial(w)<r\mbox{  if
}\sigma(w)=0,\quad\quad \partial(w)=r\mbox{ if
}\sigma(w)=-1.$$Extend $\partial$ to a function
$\partial:W\to[-r,r]\cap{\mathbb Z}$ by setting
$\partial(ws_d)=-\partial(w)$ for $w\in W^{s_d}$. Then, as we assume $d\le 2$, properties (\ref{homebed1}) and
(\ref{homebed2}) are empty resp. fulfilled for trivial reasons. Therefore we
conclude with Proposition \ref{silvester}.\hfill$\Box$\\

\begin{flushleft} \textsc{Humboldt-Universit\"at zu Berlin\\Institut f\"ur Mathematik\\Rudower Chaussee 25\\12489 Berlin, Germany}\\ \textit{E-mail address}:
gkloenne@math.hu-berlin.de \end{flushleft} \end{document}